\documentclass[reqno]{amsart}

\usepackage[english]{babel} 

\usepackage{amsmath}
\usepackage{amsfonts}
\usepackage{amssymb}
\usepackage{amsthm}
\usepackage{dsfont}

\usepackage{enumitem}

\input{xy}
\xyoption{all}

\newcommand{\N}{\mathbb{N}}
\newcommand{\R}{\mathbb{R}}
\newcommand{\zero}{{\bf 0}}

\newcommand{\one}{\mathds 1}
\newcommand{\pol}[1]{|#1|}
\renewcommand{\P}{\nabla(P)}

\newtheorem*{thm_main}{Theorem}
\newtheorem{cl}{Lemma}[section]
\newtheorem{cor}{Corollary}[section]
\newtheorem{prop}{Proposition}[section]
\newtheorem{lemma}{Claim}[section]

\theoremstyle{definition}
\newtheorem{defn}{Definition}[section]

\newtheorem{oss}{Remark}

\newcommand{\remove}[1]{}

\title[Euler characteristic and vector lattices]{The Euler characteristic of a polyhedron as a valuation on its coordinate vector lattice}

\author{Andrea Pedrini}

\address[A. Pedrini]{Dipartimento di Matematica ``F. Enriques''\\
Universit\`{a} degli Studi di Milano \\
via Cesare Saldini 50 \\ I-20133 Milano \\
Italy}

\email{andrea.pedrini@unimi.it}


\begin{document}

\nocite{*} 

\begin{abstract}
 A celebrated theorem of Hadwiger states that the Euler-Poincar\'e characteristic is the the unique invariant and continuous valuation on the distributive lattice of compact polyhedra in $\R^n$ that assigns  value one to each convex non-empty such polyhedron. This paper provides an analogue of Hadwiger's result for finitely presented unital vector lattices (i.e.\ real vector spaces with a compatible lattice order, also known as Riesz spaces). 
 
 The vector lattice of continuous and piecewise (affine) linear real-valued functions on a  compact polyhedron, with operations defined pointwise from the vector lattice $\R$,  is a finitely presented unital vector lattice; and it is a non-trivial fact that all such vector lattices arise in this manner, to within an isomorphism. Each function in such a  vector lattice can be written as a linear combination of a subset of distinguished elements that we call vl-Schauder hats. We prove here that the functional that assigns to each non-negative piecewise linear function on the polyhedron the Euler-Poincar\'e characteristic of its support is the unique vl-valuation (a special class of valuations on vector lattices) that assigns one to each vl-Schauder hat of the vector lattice.
\end{abstract}
\maketitle

\section{Introduction}\label{s:intro}
The (always compact, but not necessarily convex) polyhedra lying in the real vector space $\R^n$, for $n\geq 1$ an integer, form a distributive lattice $\mathcal{P}$ under set-theoretic union and intersection. The convex polyhedra in  $\R^n$, also known as polytopes, generate the whole lattice $\mathcal{P}$ via finite unions. A celebrated theorem of Hadwiger determines the Euler-Poincar\'e characteristic $\chi(P)$ of an element  $P\in \mathcal{P}$ as the value at $P$ of the unique (lattice-theoretic, real-valued) invariant and continuos valuation $\chi\colon\mathcal{P}\to \R$ that assigns $1$ to each non-empty polytope in $\R^n$, and $0$ to the empty polytope; for a proof, see \cite[Theorem 5.2.1]{klain_rota}. Recall that, abstractly, a  valuation on a lattice $L$ is a function $\nu \colon L \to \R$ such that $\nu(x)+\nu(y)=\nu(x\vee y)+\nu(x\wedge y)$ holds for all $x,y \in L$; see \cite[page 74]{birkhoff}. Further recall that the number $\chi(P)$ is a homotopy-invariant of  the polyhedron $P$ that can be defined as follows. Choose any triangulation $\Delta$ of  $P$, and let $\alpha_m$ be the number of faces of $\Delta$ that   have dimension $m$, for $m\geq 0$ an integer. Then
\[\tag{*}\label{t:chi}
  \chi(P) = \sum_{m=0}^n(-1)^m\alpha_m
\]
It is a classical result that the sum on the right-hand side of (\ref{t:chi}) is independent of the choice of $\Delta$. For background, see e.g.\ \cite{maunder}.

\smallskip Our purpose in this paper is to provide an analogue of Hadwiger's theorem using the unital vector lattice of (always continuous) piecewise (affine) linear functions $P\to \R$.
In the remaining part of this introduction we state the main result; all definitions not explicilty given here can be found in Section \ref{s:bg}.

\smallskip 
We assume familiarity with  vector lattices, also known as \emph{Riesz spaces}; standard references are  \cite{birkhoff} and \cite{luxembourg_zaanen}.
We recall that a \emph{\textup{(}real\textup{)} vector lattice}  is an algebra ${\bf V}=(V,+,\land,\vee,\{\lambda\}_{\lambda\in\R},\zero)$ such that
 $(V,+,\{\lambda\}_{\lambda\in\R},\zero)$ is a (real) vector space,
$(V,\land,\vee)$ is a lattice, addition distributes over meets (i.e.\ for all $t,v,w\in V$ we have $t+(v\land w)=(t + v) \land (t + w)$), and so does multiplication by positive scalars (i.e.\
  for all $v,w\in V$ and for all $\lambda \in \R$, if $\lambda\geq0$ then $\lambda(v\land w)=\lambda v\land\lambda w$). Thus vector lattices form a variety of algebras (with continuum-many operations) by their very definition.
  It is well known that the underlying lattice of ${\bf V}$ is necessarily distributive. The \emph{positive cone} of ${\bf V}$ is the collection ${\bf V}^+=\{v \in V \mid v\geq \zero\}$ of its non-negative elements, where $\geq$ denotes the underlying lattice order of ${\mathbf V}$. Morphisms of vector lattices are homomorphisms in the variety, that is, linear maps that also preserve the lattice structure.
 From now on we shall follow common practice and blur the distinction between {\bf V} and its underlying set $V$. Moreover, we will denote both the element $\zero$  of $V$ and the real number zero by the same symbol $0$: the meaning will be clear from the context.
 
 \smallskip We can now define the new notion introduced in this paper.
\begin{defn}[Vl-valuation]
Let $V$ be a  vector lattice with underlying set $V$,  and let  $V^+$ be its positive cone.
 A \emph{vl-valuation} on $V$ is a function $\nu\colon V \to\R$ such that:
 \begin{enumerate}[label=\emph{V\arabic*)}, leftmargin=1.5cm]
  \item $\nu(0)=0$,
  \item for all $x,y \in V$, $\nu(x)+\nu(y)=\nu(x \vee y)+\nu(x\land y)$,
  \item for all $x,y \in V^+$, $\nu(x+ y)=\nu(x \vee y)$,
  \item for all $x,y \in V^+$, if $x\land y=0$ then $\nu(x-y)=\nu(x)-\nu(y)$.
 \end{enumerate}
\end{defn}
 \noindent We shall be concerned with vl-valuations on a 
 class of vector lattices naturally arising from polyhedra, as follows. 
Given a set $S$ in $\R^n$, a function $f:S\to\R$ is said to be \emph{piecewise linear} if there is a finite set $f_1,\dots,f_m$ of affine linear functions $\R^n\to\R$ such that for all $s\in S$ there exists an index $i\leq m$ for which $f(s)=f_i(s)$. 
 Now we consider a  polyhedron $P$ in $\R^n$ and the set of all functions $f:P\to\R$ that are continuous with respect to the Euclidean metric and piecewise linear, equipped with pointwise addition, supremum, infimum, products by real scalars and the zero function. It is easy to show that this set is a vector lattice under the mentioned operations;  we will denote it  $\nabla{(P)}$.

\smallskip 
Given a vector lattice $V$, an element $u\in V$ is a \emph{strong order unit}, or just a \emph{unit} for short, if for all $0\leq v\in V$ there exists a $0\leq\lambda\in\R$ such that $v\leq\lambda u$.  A \emph{unital vector lattice} is a pair $(V,u)$, where $V$ is a vector lattice and $u$ is a unit of $V$. It turns out that every finitely generated vector lattice admits a unit. For if $v_1,\ldots,v_u$ is a finite set of generators for $V$, then it is easily checked that $|v_1|+\cdots+|v_u|$ is a unit of $V$, where $|v_i|=v_i\vee -v_i$ as usual. Morphisms of unital vector lattices are the vector-lattice homomorphisms that carry units to units. Such homomorphisms are called  \emph{unital}. Note that unital vector lattices do not form a variety of algebras, because the Archimedean property of the unit is not even definable by first-order formul\ae, as is shown via a standard compactness argument.
\smallskip

It is clear that the  element $\one \in \nabla(P)$ defined as the function identically equal to $1$ on the polyhedron $P$ is  a  unit of $\nabla(P)$. Hence the pair $(\nabla(P),\one)$ is  a unital vector lattice. Moreover, the unital vector lattices of the form $(\nabla(P),\one)$ can be characterized abstractly. Recall that a vector lattice is finitely presented if it is the quotient of the free vector lattice on $n$ generators $\text{FVL}(n)$ by a finitely generated ideal.
(Free vector lattices exist because, as we observed above, vector lattices form a variety of algebras.) The well-known \emph{Baker-Beynon duality} states that the category of finitely presented vector lattices with vector lattice morphisms is dually equivalent to the category of polyhedral cones (with vertex at $0$) in some Euclidean space and with piecewise homogeneous linear continuous maps. Moreover, there is an induced duality between the category of finitely presented vector lattices with a distinguished strong unit and unital morphisms and the category of polyhedra and piecewise linear continuous maps, see \cite{baker} and \cite{beynon} for more details. This  duality entails that finitely presented unital vector lattices are exactly those representable as $(\nabla(P),\one)$ to within a unital isomorphism, for some polyhedron $P$ in some Euclidean space $\R^n$. 
The polyhedron $P$ is called the \emph{support} of the vector lattice itself, and $(\nabla(P),\one)$ is its \emph{coordinate vector lattice}.

\smallskip
We next isolate a special class of elements of $(\nabla(P),\one)$ that we call \emph{vl-Schauder hats} (or just \emph{hats}, for short). Hats
form a generating set  of the underlying vector space of the vector lattice: each element of $\nabla(P)$ is a linear combination of a suitable finite set of hats. In the present  context, vl-Schauder hats  play the same role as Schauder hats in lattice-ordered Abelian groups \cite[and references therein]{manara_marra_mundici} and MV-algebras  \cite[9.2.1]{cignoli_dottaviano_mundici}. Pursuing the analogy with Hadwiger's above-mentioned theorem,  we will see in due course that our version of the Euler-Poincar\'e characteristic assigns value one to each hat.

Formally, we define as follows (please see Section~\ref{s:bg} for background on triangulations).
\begin{defn}[vl-Schauder hats]
\label{def:hats}
A \emph{vl-Schauder hat} is an element $h\in(\P,\one)$ for which there are a triangulation $K_h$ of $P$ linearizing $h$ and a vertex $\bar{x}$ of $K_h$ such that $h(\bar{x})= 1$ and $h(x)=0$ for any other vertex $x$ of $K_h$. 
\end{defn}

We remark that it is possible to characterize vl-Schauder hats abstractly in the language of vector lattices, transposing to our context the results  in \cite{marra2}. Therefore our main result, stated below, may be regarded as a theorem about unital vector lattices that does not depend on any geometrical representation.

%
\begin{thm_main}\label{t:main}Let $P$ be a polyhedron in $\R^n$, for some integer $n\geq 1$, and let $(\P,\one)$ be the
finitely presented unital vector lattice  of real-valued piecewise linear functions on $P$. Then there is a unique vl-valuation 
\[
\alpha\colon \P \to \R
\]
that assigns the value $1$ to each vl-Schauder hat in $\P$. Moreover,  the number $\alpha(\one)$ is the Euler-Poincar\'e characteristic \textup{(\ref{t:chi})}
of the polyhedron $P$.
\end{thm_main}
In fact, our proof of Theorem \ref{t:main} in Section \ref{s:proofmain} below will show more: for each non-negative element $f \in \P$, $\alpha(f)$ is the Euler-Poincar\'e characteristic of the support of $f$. (By the support of $f\colon P \to \R$ we mean the complement of the vanishing locus of $f$, as usual). Since the support is an open set, it is not in general compact and therefore cannot be triangulated by a finite simplicial complex. Thus the classical combinatorial formula \textup{(\ref{t:chi})} cannot be used to define the characteristic of the support of $f$. Nonetheless, we can use the \emph{supplement} of the support of $f$. The supplement (see Section~\ref{s:bg}) is a standard construction in algebraic topology: it is a simplicial complex $\bar{L}$ that approximates the set-theoretic difference between the underlying polyhedra $|K|$ and $|L|$ of a simplicial complex $K$ and its subcomplex $L$. It can be shown (see~\cite[Proposition 5.3.9]{maunder}) that $|\bar{L}|$ is homotopically equivalent to the set-theoretic difference $|K|\backslash|L|$. So the Euler-Poicar\'e characteristic of  $|\bar{L}|$ given by \textup{(\ref{t:chi})} is exactly the Euler-Poincar\'e characteristic of $|K|\backslash|L|$, defined in a standard way as the homological generalization of \textup{(\ref{t:chi})} through the formula
\begin{equation*}
  \chi(|K|\backslash|L|) = \sum_{m=0}^\infty(-1)^m\beta_m,
\end{equation*}
where each $\beta_m$ is the rank of the $m$th homology group $H_m(|K|\backslash|L|)$.

\section{Preliminaries and background}\label{s:bg}

\subsection{Simplicial complexes and polyhedra}

For more details on polyhedra and proofs of the results quoted here, see \cite{maunder}.

\smallskip
Given $n+1$ affinely independent points $x_0,\dots,x_n$ in some Euclidean space $\R^m$, we say that the \emph{$n$-simplex} $\sigma_n=(x_0,\dots,x_n)$ is the set of all the convex combinations of $x_0,\dots,x_n$, that is the set of points $\sum_{i=0}^n\lambda_i x_i$, where the $\lambda_i$ are real numbers such that $\lambda_i\geq0$ for all $i$ and $\sum_{i=0}^n\lambda_i=1$. The points $x_0,\dots,x_n$ are called the \emph{vertices} of $\sigma_n$, and the number $n$ is the \emph{dimension} of $\sigma_n$.
  The (\emph{relative}) \emph{interior} of $\sigma_n$ is the subset of $\sigma_n$ of those points $\sum_{i=0}^n\lambda_i x_i$ such that $\lambda_i>0$ for all $i$. The \emph{barycentre} of $\sigma_n$ is the point
 \begin{equation*}
  \hat{\sigma}_n = \left(\frac{1}{n+1}\right)(x_0+\dots+x_n).
 \end{equation*}
 
 A \emph{face} of $\sigma_n$ is the convex hull of a subset of its vertices. The \emph{boundary} $\dot{\sigma}_n$ of $\sigma_n$ is the set of all faces of $\sigma_n$ other than $\sigma_n$ itself.

\smallskip
A \emph{simplicial complex} $K$ is a finite set of simplices such that
\begin{enumerate}
 \item if $\sigma_n\in K$ and $\tau_p$ is a face of $\sigma_n$, then $\tau_p\in K$,
 \item if $\sigma_n,\tau_p\in K$, then $\sigma_n\cap\tau_p$ is either empty or is a common face of $\sigma_n$ and $\tau_p$.
\end{enumerate}
The \emph{dimension} of $K$ is the maximum dimension of its simplices. 
A \emph{subcomplex} $L$ of $K$ is a subset of simplices of $K$ that is itself a simplicial complex.
For each $r\geq0$, the \emph{$r$-skeleton} $K^r$ of $K$ is the subset of its simplices of dimension at most $r$: clearly it is a subcomplex of $K$.
The \emph{polyhedron} $\pol{K}$ of $K$ is the set of points of $\R^m$ that lie in at least one of the simplices of $K$, topologized as a subspace of $\R^m$.

\begin{defn}[Polyhedron]
 A \emph{polyhedron} is a subset $P$ of $\R^m$ such that there exists a simplicial complex $K$ with $\pol{K}=P$. 
 \end{defn}
 Each simplicial complex $K$ with polyhedron $P$ is called a \emph{triangulation} of $P$. It can be proved that any two finite triangulations of $P$ have the same dimension. Hence the \emph{dimension} of $P$ is defined to be the dimension of any one of its triangulations.
  Given a triangulation $K$ of the polyhedron $P$, a \emph{refinement} of $K$ is a triangulation $K^*$ of $P$ such that for all simplices $\sigma$ of $K^*$ there is a simplex $\tau$ of $K$ such that $\pol{\sigma}\subseteq\pol{\tau}$. Moreover, given any two triangulations $K_1$ and $K_2$ of the same polyhedron $P$, there always exists a common refinement $K^*$ of $K_1$ and $K_2$.

\begin{prop}[{\cite[Proposition 2.3.6]{maunder}}]
\label{prop:interior}
 If $K$ is a simplicial complex, then each point $x$ of $\pol{K}$ is in the interior of exactly one simplex $\sigma^x$ of $K$.
\end{prop}
The simplex $\sigma^x$ of Proposition~\ref{prop:interior} is the inclusion-smallest simplex containing $x$, and it is called the \emph{carrier} of $x$.

\smallskip

Given a  simplicial complex $K$, its \emph{derived complex} $K'$ is its first barycentric subdivision, see \cite[p.\ 48--49]{maunder}.
If $L$ is a subcomplex of $K$, it is clear that $L'$ is a subcomplex of $K'$.
The \emph{supplement} $\bar{L}$ of $L$ in $K$ is the subcomplex of $K'$ consisting of those simplices having no vertex in $L'$.

%

\subsection{Vector lattices and triangulations}

 The \emph{positive cone} of a vector lattice $V$ is the set
 $$V^+=\{x\in V : x\geq 0\}.$$
 
 Moreover, for all $x\in V$ there are two elements $x^+=x\vee0$, $x^-=(-x)\vee0$ in the positive cone $V^+$ such that $x=x^+-x^-$, and $x^+\land x^-=0$. Here, $x^+$ is the \emph{positive part} of $x$, $x^-$ is its \emph{negative part}. If $x=a-b$, with $a,b\in V^+$ and $a\land b=0$, then $a=x^+$ and $b=x^-$.

Fix a polyhedron $P$.
 Given a function $f\in\P$, we say that a \emph{linearizing triangulation for $f$} is a triangulation $K_f$ of $P$ such that $f$ is linear on $\pol{\sigma}$, for all simplices $\sigma\in K_f$.
 
 \begin{oss}
 \label{oss:linfg}
Because $f$ is piecewise linear, a standard argument shows that there always exists a linearizing triangulation $K_f$ of $P$. Indeed, if $f,g\in\nabla(P)$, there exists a triangulation of $P$ that is linearizing for both $f$ and $g$.  Moreover, if $K^*$ is a refinement of a linearizing triangulation $K$ for $f$, then $K^*$ is linearizing for $f$. Henceforth, $K_f$ will denote a linearizing triangulation of $P$ for $f$.
\end{oss}


We observe the following.

\begin{lemma}
Let $f\in\P^+$ and let $Z_{K_f,f}$ be the subcomplex of $K_f$ of those simplices where $f$ is identically zero. Then $\pol{Z_{K_f,f}}$ is the zero-set $f^{-1}(0)$ of $f$, and so does not depend on the particular choice of the triangulation $K_f$ that linearizes $f$.  
\begin{proof}
We have to prove that, given a linearizing triangulation $K_f$, $\pol{Z_{K_f,f}}$ is the zero set of $f$, that is the set $\{y\in P : f(x)=0\}$.
The inclusion of $\pol{Z_{K_f,f}}$ in the zero set of $f$ is trivial: if $y\in\pol{Z_{K_f,f}}$, then $y$ is a point of $P$ that lies in at least a simplex of $Z_{K_f,f}$ and so $f(y)=0$.
For the inverse inclusion, we have that if $y$ is a point of $P=\pol{K_f}$ such that $f(y)=0$, then, by Proposition \ref{prop:interior}, there is a simplex $\sigma$ of $K_f$ such that $y$ is a point of the interior of $\sigma$. Recalling that $f\geq0$, the linearity of $f$ on the simplices of $K_f$, and in particular on $\sigma$, ensures that $f$ is identically zero on the whole simplex $\sigma$. So $\sigma$ is a simplex of $Z_{K_f,f}$ and $y\in\pol{Z_{K_f,f}}$.
\end{proof}
\end{lemma}

Defining vl-Schauder hats as in Definition~\ref{def:hats}, given a triangulation $K$ of $P$ with vertices $\{x_0,\dots,x_m\}$, the \emph{vl-Schauder hats of $K$} are those vl-Schauder hats $\{h_i\}$ such that, for all $i,j\in\{0,\dots,m\}$, $h_i(x_i)=1$ and $h_i(x_j)=0$. The uniquely determined $x_i$ such that $h_i(x_i)=1$ is called the \emph{vertex} of $h_i$.  Each $f\in\P$ can be written as a sum $\sum_{i=0}^m a_ih_i$ (where $a_i\in\R$) of distinct vl-Schauder hats $h_0,\dots,h_m$ of a common linearizing triangulation $K_f$ for $f$. If $f\in\P^+$, then necessarily $0\leq a_i$ for all $i=0,\dots, m$.

\begin{oss}
\label{oss:vl-Schauderhats}
Given any two vl-Schauder hats $h_i$ and $h_j$ of a triangulation $H$, with vertices $x_i$ and $x_j$, respectivly, the element $h_{ij}=2(h_i \land h_j)$ is either zero or a vl-Schauder hat with vertex $x_{ij}=\widehat{(x_i,x_j)}$. Hence $h_{ij}$ is a vl-Schauder hat of the triangulation $H^*$ obtained from $H$ by performing the barycentric subdivision of the simplex $(x_i,x_j)$. Explicity, by adding the $0$-simplex $(x_{ij})$ to $H$, and replacing each $n$-simplex of the form $\sigma=(x_{k_0},\dots,x_i,\dots,x_j,\dots,x_{k_n})$ with the $n$-simplices $\tau=(x_{k_0},\dots,x_i,\dots,x_{ij},\dots,x_{k_n})$ and $\rho=(x_{k_0},\dots,x_{ij},\dots,x_j,\dots,x_{k_n})$. The vl-Schauder hats associated with $x_i$ and $x_j$ in the new triangulation $H^*$ are, respectively, $h_i'=h_i-(h_i \land h_j)$ and $h_j'=h_j-(h_i \land h_j)$.
\end{oss}

\subsection{Pc-valuations}
%

We presently show that a vl-valuation is uniquely determined by its values at the positive cone.

\begin{defn}[Pc-valuation]
 Let $V$ by a vector lattice.
 A \emph{pc-valuation} on the positive cone is a function $\nu^+:V^+\to\R$ such that:
 \begin{enumerate}[label=\emph{P\arabic*)}, leftmargin=1.5cm]
  \item $\nu^+(0)=0$,
  \item for all $x,y \in V^+$, $\nu^+(x \vee y)=\nu^+(x)+\nu^+(y)-\nu^+(x\land y)$,
  \item for all $x,y \in V^+$, $\nu^+(x+ y)=\nu^+(x \vee y)$.
 \end{enumerate}
\end{defn}

\begin{cl}
\label{pcval}
The operation of restriction of a vl-valuation to the positive cone is a bijection between the set of all vl-valuations on $V$ and the set of all pc-valuations on $V^+$. The inverse bijection is the operation that extends a pc-valuation $\nu^+$ to the vl-valuation 
\begin{equation}
\label{eq:pcext}
 \nu^\pm: x \mapsto \nu^+(x^+)-\nu^+(x^-),
\end{equation} 
where $x^+$ and $x^-$ are, respectively, the positive and the negative part of $x$.
\begin{proof}
Clearly, if $\nu$ is a vl-valuation on $V$, then its restriction $\nu|_{V^+}$ is a pc-valuation. On the other hand, if we consider a pc-valuation $\nu^+$ defined on $V^+$, then its extension $\nu^\pm$ given in~\eqref{eq:pcext} is a vl-valuation:
\begin{enumerate}[label=V\arabic*), leftmargin=1.5cm]
  \item $\nu^\pm(0)=\nu^+(0^+)-\nu^+(0^-)=\nu^+(0)-\nu^+(0)=0$;
  \item for all $x,y \in V$, we have $(x\vee y)^+=(x^+\vee y^+)$, $(x\land y)^+=(x^+\land y^+)$, $(x\vee y)^-=(x^-\land y^-)$, $(x\land y)^-=(x^-\vee y^-)$ and so 
  \begin{equation*}
  					\begin{split}
					\nu^\pm(x \vee y)&=\nu^+((x \vee y)^+)-\nu^+((x \vee y)^-)=\\
					&=\nu^+(x^+ \vee y^+)-\nu^+(x^- \land y^-) =\\
					&=\nu^+(x^+)+\nu^+(y^+)-\nu^+(x^+\land y^+)\\
					  &\qquad+\nu^+(x^-\vee y^-)-\nu^+(x^-)-\nu^+(y^-)=\\
					&=\nu^+(x^+)-\nu^+(x^-)+\nu^+(y^+)-\nu^+(y^-)\\
					  &\qquad-(\nu^+((x\land y)^+)-\nu^+((x\land y)^-))=\\
					&=\nu^\pm(x)+\nu^\pm(y)-\nu^\pm(x\land y);
					\end{split}
				        \end{equation*}
  \item for all $x,y \in V^+$, we have $x+y \in V^+$, $(x+y)^+=x+y$ and $(x+y)^-=0$, and so \begin{equation*}
  					\begin{split}
					\nu^\pm(x+ y)&=\nu^+((x+ y)^+)-\nu^+((x+ y)^-)=\\
					&=\nu^+(x+ y)-\nu^+(0)=\\
					&=\nu^+(x\vee y)-0=\\
					&=\nu^+(x^+\vee y^+)=\\
					&=\nu^+((x\vee y)^+)-0=\\
					&=\nu^+((x\vee y)^+)-\nu^+((x\vee y)^-)=\\
					&=\nu^\pm(x\vee y);
					\end{split}
				        \end{equation*}
  \item for all $x,y \in V^+$, if $x\land y=0$, then we have $x^-=y^-=0$, $(x-y)^+=x^+$ and $(x-y)^-=y^+$, and so \begin{equation*}
  					\begin{split}
					\nu^\pm(x-y)&=\nu^+((x-y)^+)-\nu^+((x-y)^-)=\\
					&=\nu^+(x^+)-\nu^+(y^+)=\\
					&=\nu^+(x^+)-\nu^+(x^-)-(\nu^+(y^+)-\nu^+(y^-))=\\
					&=\nu^\pm(x)-\nu^\pm(y).
					\end{split}
				        \end{equation*}
 \end{enumerate}
 Moreover, if $\nu$ is a vl-valuation, then, for all $x\in V$,
 \begin{equation*}
 \begin{split}
   (\nu|_{V^+})^\pm(x)&=\nu|_{V^+}(x^+)-\nu|_{V^+}(x^-)=\\
   &=\nu(x^+)-\nu(x^-)=\\
   &=\nu(x^+-x^-)=\nu(x).
   \end{split}
 \end{equation*}
 On the other hand, if $\nu^+$ is a pc-valuation with extension $\nu^\pm$, then, for all $x\in V^+$,
 \begin{equation*}
   \nu^\pm|_{V^+}(x)=\nu^\pm(x)=\nu^+(x^+)-\nu^+(x^-)=\nu^+(x^+)=\nu^+(x).
 \end{equation*}
\end{proof}
\end{cl}

So, from now on, without loss of generality, we can consider only pc-valuations and positive cones.

\begin{cl}
\label{cl:vhats}
 If $\nu^+$ is a pc-valuation and $x,y\in V^+$, then, for all $0<a\in\R$,
 $$\nu^+(x + ay) = \nu^+(x + y).$$
 \begin{proof}
 Let $0<b\in\R$, then, for all $\dfrac{b}{2}\leq c\leq b$, $0\leq b-c\leq\dfrac{b}{2}$ and so
  \begin{equation*}
  \begin{split}
  \nu^+(x+by)&=\nu^+((x +cy)+(b-c)y)=\\
  &=\nu^+(x+cy)+\nu^+((b-c)y)-\nu^+((x+cy)\land(b-c)y)=\\
  &=\nu^+(x+cy)+\nu^+((b-c)y)-\nu^+((b-c)y)=\nu^+(x+cy).
  \end{split}
  \end{equation*}
  By induction, $\nu^+(x+by)=\nu^+(x+cy)$ for all $b\geq c\geq\dfrac{b}{2^n}$ (for all $n\in\N$) and so for all $b\geq c>0$. Then we can choose $b=\max\{a,1\}$ and we have $\nu^+(x+ay)=\nu^+(x+y)$.
 \end{proof}
\end{cl}

In particular, a pc-valuation $\nu^+$ is a positively homogeneous function of degree $0$: $\nu^+(ax)=\nu^+(x)$ for all $x\in V^+$ and for all real coefficients $a>0$. More generally:

\begin{cor}
\label{cor:vf}
If $\nu^+$ is a pc-valuation on $V^+$, and if $x=\sum_{i=0}^m a_ix_i$ is such that  $0<a_i\in\R$ and $x_0,\dots x_m\in V^+$, then
\begin{equation*}
 \nu^+(x)=\nu^+\left(\sum_{i=0}^m x_i\right).
\end{equation*}
\end{cor}

\section{Proof of Theorem}\label{s:proofmain}
To prove Theorem~\ref{s:proofmain}, we will caracterize exactly the vl-valuation that assigns value one to each vl-Schauder hat. In light of Lemma~\ref{pcval}, we can restrict attention to pc-valuation on the positive cone.

First of all we observe that, as suggested by Lemma~\ref{cl:vhats}, a pc-valuation forgets the height of functions, so the only information it retains is concerned with supports and vanishing loci. We therefore try to use the Euler-Poincar\'e characteristic of the support $\text{supp}(f)$ of the functions in $\P^+$ to construct our pc-valuation. To this aim, we use the supplement. For each $f\in\P^+$, we choose a linearizing triangulation $K_f$ and isolate the zero-set $Z_{K_f,f}$ of $f$. Then we build its supplement $\overline{Z_{K_f,f}}$ in $K_f$. Now we have a simplicial complex (and so an associated polyhedron) which approximates the support of $f$, and we can compute its Euler-Poincar\'e characteristic. Recalling that $|Z_{K_f,f}|$ is the vanishing locus of $f$, and that by~\cite[Proposition 5.3.9]{maunder} there is a homotopy equivalence between $|\overline{Z_{K_f,f}}|$ and $\text{supp}(f)=P\backslash f^{-1}(0)=|K_f|\backslash|Z_{K_f,f}|$, we have that the Euler-Poicar\'e characteristic of $|\overline{Z_{K_f,f}}|$ does not depend on the choice of the linearizing triangulation $K_f$, but just on the homotopy type of $|\overline{Z_{K_f,f}}|$, and so only on the homotopy type of $\text{supp}(f)$. The Euler-Poincar\'e characteristic, as defined in \textup{(\ref{t:chi})}, indeed is a homotopy invariant (see~\cite[Lemma 4.5.17]{maunder} and remarks following it).

Because of this, the following is well defined.

\begin{defn}
We define $\alpha^+:\P^+\to\R$ as
$$\alpha^+(f)=\chi(\text{supp}(f))=\chi(|\overline{Z_{K_f,f}}|),$$
where $K_f$ is a linearizing triangulation for $f\in\P^+$, and $\chi$ is the Euler-Poincar\'e characteristic defined in \textup{(\ref{t:chi})}.
\end{defn}

We now prove that $\alpha^+$ is a pc-valuation that assigns $1$ to each vl-Schauder hat.

\begin{cl}
\label{cl:alfavl-Schauder}
The following hold.
\begin{enumerate}
 \item $\alpha^+(0)=0$;
 \item if $h\in \P^+$ is a vl-Schauder hat, then $\alpha^+(h)=1$;
 \item for all $f,g\in \P^+$,
 $$\alpha^+(f + g)=\alpha^+(f \vee g)=\alpha^+(f)+\alpha^+(g)-\alpha^+(f\land g).$$
\end{enumerate}
\begin{proof}
 (1)\  $\alpha^+(0)=\chi(\text{supp}(0))=\chi(\emptyset)=0$.
 
 \smallskip 
 \noindent
 (2)\  We can choose a triangulation $K_h$ that linearizes $h$ and such that the vertex $\tilde{x}$ is the only one on which $h>0$. Then, we observe that $\overline{Z_{K_h,h}}$ is the \emph{simplicial neighbourhood} $N_{K_h'}(\tilde{x})$ of $\tilde{x}$ in $K_h'$, that is the smallest subcomplex of $K_h'$ containing each simplex of $K_h'$ which contains $\tilde{x}$. It can be shown (using, for example, \cite[Proposition 2.4.4]{maunder}) that $|N_{K_h'}(\tilde{x})|$ is contractible, that is it is homotopically equivalent to the point $\tilde{x}$: so the Euler-Poincar\'e characteristic of $|\overline{Z_{K_h,h}}|$ is the same of the single point $\tilde{x}$. This proves that $\alpha^+(h)=1$.
 
 \smallskip 
 \noindent
 (3)\  By Remark~\ref{oss:linfg}, we can always chose a triangulation $K$ of $P$ that simultaneously linearizes $f + g$, $f \vee g$, $f$, $g$ and $f \land g$. Let us compute the Euler-Poincar\'e characteristic using this common linearizing triangulation.
 
 First, we observe that, by the linearity of $f$ on $K$, on the barycentres of the simplices of $Z_{K,f}$ the function $f$ takes value $0$. Hence $f$ is identically $0$ on each vertex of $Z_{K,f}'$. Furthermore, again by linearity, the vertices of $Z_{K,f}'$ are exactly all vertices of $K'$ where $f$ is $0$. So, if we compute $\overline{Z_{K,f}}$ as the set of simplices of $K'$ with no vertices in $Z_{K,f}'$, we have that it is the subset of $K'$ of all those simplices whose vertices are in the support of $f$, and so, by linearity,
 \begin{equation*}
  \overline{Z_{K,f}}=\{\sigma\in K'\mid \sigma\subseteq\text{supp}(f)\},
 \end{equation*}
 and similarly for $g$.
 The same reasoning can be applied to $f\land g$, $f\vee g$ and $f+g$. Observing that $\text{supp}(f\land g)=\text{supp}(f)\cap\text{supp}(g)$ and $\text{supp}(f+g)=\text{supp}(f\vee g)=\text{supp}(f)\cup\text{supp}(g)$, we have:
 \begin{equation*}
  \overline{Z_{K,f\land g}}=\{\sigma\in K'\mid \sigma\subseteq\text{supp}(f)\cap\text{supp}(g)\},
 \end{equation*}
 \begin{equation*}
  \overline{Z_{K,f\vee g}}=\{\sigma\in K'\mid \sigma\subseteq\text{supp}(f)\cup\text{supp}(g)\},
 \end{equation*}
  \begin{equation*}
  \overline{Z_{K,f+g}}=\{\sigma\in K'\mid \sigma\subseteq\text{supp}(f)\cup\text{supp}(g)\}.
 \end{equation*}
 Because $\overline{Z_{K,f\vee g}}=\overline{Z_{K,f+ g}}$ the first equality $\alpha^+(f + g)=\alpha^+(f \vee g)$ is trivial. For the second one, we have that the $m$-simplex $\sigma_m$ is in $\overline{Z_{K,f\land g}}$ if, and only if, it is in both $\overline{Z_{K,f}}$ and $\overline{Z_{K,g}}$, and, in this case, it also lies in $\overline{Z_{K,f\vee g}}$. So, if $\alpha_{m,\star}$ is the number of $m$-simplices in $\overline{Z_{K,\star}}$, we have $\alpha_{m,f \vee g}=\alpha_{m,f}+\alpha_{m,g}-\alpha_{m,f\land g}$. Summing over $m$ completes the proof.
\end{proof}
\end{cl}

\begin{oss}
Observe that $\alpha(\one)=\chi(P)$. In fact, each triangulation $K$ of $P$ is a linearizing triangulation for $\one$ and $Z_{K,\one}=\emptyset$; then $\overline{Z_{K,\one}}=K'$ and so $\pol{\overline{Z_{K,\one}}}=\pol{K'}=P$.
\end{oss}

The following technical result is crucial.

\begin{cl}
\label{cl:andhat}
 Let $h_0,\dots,h_n$ be distinct vl-Schauder hats of the same triangulation $H$ of $P$, then
 \begin{equation}
 \label{eq:andhat}
  \left(\sum_{i=0}^{n-1}h_i\right)\land h_n=\sum_{i=0}^{n-1}k_i,
 \end{equation}
 where
 \begin{align*}
  &k_i=h_i\land h_n^i,\\
  &h_n^0=h_n,\\
  &h_n^i=h_n^{i-1}-(h_{i-1} \land h_n^{i-1})=h_n^{i-1}- k_{i-1}.
 \end{align*}
 Moreover, the non-zero elements of the set $\{2k_0,\dots,2k_{n-1}\}$ are distinct vl-Schauder hats of a single triangulation $K$ of $P$.
 
 \begin{proof}
  First we notice that:
  \begin{itemize}
   \item $h_n^i\leq h_n$,
   \item $k_i\leq h_i$ and $k_i\leq h_n$,
   \item if $h_n^i=0$, then $\forall j\geq i$ $h_n^j=0$ and $k_j=0$,
   \item if $2k_i\neq0$ and $2k_j\neq0$ (with $i\neq j$), then they are distinct: in fact, recalling  Remark~\ref{oss:vl-Schauderhats}, $h_n^l$ is always a vl-Schauder hat associated with the point $x_n$, and so $2k_i$ attains its maximum at $x_{in}$, but the maximum of $2k_j$ is attained at $x_{jn}$ and $x_{in}\neq x_{jn}$ because $h_i\neq h_j$.
  \end{itemize}
  The proof proceeds by induction on $n$. 
   If $n=1$, there is nothing to prove. The only thing we need to observe is that $2k_0=2(h_0\land h_1)$ is either zero or a vl-Schauder hat of the triangulation $H^*=K$ given in Remark~\ref{oss:vl-Schauderhats}. 
   
  Assume the thesis be true for all $m<n$. In particular, $$\left(\sum_{i=0}^{n-2}h_i\right)\land h_n=\sum_{i=0}^{n-2}k_i,$$
   where $2k_0,\dots,2k_{n-2}$ are either the zero function or vl-Schauder hats of the single triangulation $\tilde{K}$, together with the hats $h_i'=h_i-k_i$ for $i=0,\dots,n-2$ and $h_n^{n-1}$.
   Thanks to Remark~\ref{oss:linfg}, we can now take a new triangulation $L$ that linearizes all the functions involved in the proof and then consider the restrictions of these functions on each single simplex of $L$.
   There are three cases:

\smallskip
\noindent
    \emph{Case 1.}\  $h_n\leq \sum_{i=0}^{n-2}h_i\leq\sum_{i=0}^{n-1}h_i$. In this case
    $$\left(\sum_{i=0}^{n-1}h_i\right)\land h_n=h_n=\left(\sum_{i=0}^{n-2}h_i\right)\land h_n=\sum_{i=0}^{n-2}k_i,$$
    and so the only thing to prove is $k_{n-1}=0$.
     If $\exists j\in\{0,\dots,n-2\}$ such that $h_n^j\leq h_j$, then $k_j=h_n^j$ and $\forall i>j$ $h_n^i=k_i=0$; in particular $k_{n-1}=0$.
     Else, if $\forall i\in \{0,\dots, n-2\}$ $h_i<h_n^i\leq h_n$, then $k_i=h_i$ and $h_n^i=h_n-\sum_{j=1}^{i-1}h_i$. So 
     $$h_n^{n-1}=h_n-\sum_{i=0}^{n-2}h_i\leq\sum_{i=0}^{n-2}h_i-\sum_{i=0}^{n-2}h_i=0,$$
     and then $k_{n-1}=0$.

\smallskip
\noindent
\emph{Case 2.}\  $\sum_{i=0}^{n-2}h_i\leq\sum_{i=0}^{n-1}h_i<h_n$. In this case $\forall i\in\{0,\dots,n-1\}$ we have $h_i< h_n$ and $\forall j\in\{0,\dots,n-1\}$ we have $\sum_{i=0}^jh_i<h_n$. Moreover, $\forall i\in\{0,\dots,n-1\}$ we have $h_i<h_n^i$. If, indeed, there were a first index $j$ (necessarily strictly greater than $0$) such that $h_n^j\not>h_j$, we would have $h_n^j\leq h_j$, because the triangulation $H$ is supposed to be fine enough to linearize $k_j=h_n^j\land h_j$, and also $h_i<h_n^i$ for all $i<j$, and therefore $k_i=h_i$, $h_n^{i+1}=h_n-\sum_{l=0}^{i}h_l$. So $h_n^j=h_n-\sum_{i=0}^{j-1}h_i$, and then $0\leq h_n^j- h_j=h_n-\sum_{i=0}^{j}h_i$;
   eventually, we reach the contradiction $h_n\leq\sum_{i=0}^{j}h_i$. Then $k_{n-1}=h_{n-1}\land h_n^{n-1}=h_{n-1}$, and so
    \begin{equation*}
     \begin{split}
     \left(\sum_{i=0}^{n-1}h_i\right)\land h_n&=\sum_{i=0}^{n-1}h_i=\sum_{i=0}^{n-2}h_i+ h_{n-1}=\\
     &=\left(\left(\sum_{i=0}^{n-2}h_i\right)\land h_n\right)+ h_{n-1}=\\
     &=\sum_{i=0}^{n-2}k_i+ k_{n-1}=\sum_{i=0}^{n-1}k_i.
     \end{split}
    \end{equation*}

\smallskip
\noindent
\emph{Case 3.}\  $\sum_{i=0}^{n-2}h_i<h_n\leq\sum_{i=0}^{n-1}h_i$. As in the previous case, $\forall i\in\{0,\dots,n-2\}$ we have $h_i<h_n^i$ and so $h_n^i=h_n-\sum_{j=0}^{i-1}h_j$ and $k_i=h_i$. Moreover $h_n^{n-1}=h_n-\sum_{i=0}^{n-2}h_i$ and $$k_{n-1}=h_{n-1}\land\left(h_n-\sum_{i=0}^{n-2}h_i\right).$$ Therefore 
   \begin{equation*}
   \begin{split}
    \left(\sum_{i=0}^{n-1}h_i\right)\land h_n&=\left(\sum_{i=0}^{n-1}h_i\right)\land \left(\sum_{i=0}^{n-2}h_i + h_n-\sum_{i=0}^{n-2}h_i\right)=\\  
    &=\left(\sum_{i=0}^{n-2}h_i+ h_{n-1}\right)\land \left(\sum_{i=0}^{n-2}h_i + \left(h_n-\sum_{i=0}^{n-2}h_i\right)\right)=\\
    &=\sum_{i=0}^{n-2}h_i+ \left(h_{n-1}\land \left(h_n-\sum_{i=0}^{n-2}h_i\right)\right)=\\
    &=\sum_{i=0}^{n-2}k_i+ k_{n-1}=\sum_{i=0}^{n-1}k_i.
   \end{split}
   \end{equation*}
   
   \smallskip
   This proves~\eqref{eq:andhat}.
  Finally, we show that $\{2k_0,\dots,2k_{n-1}\}$ (when non-zero) are distinct vl-Schauder hats of the same triangulation. Take $\tilde{K}$, and construct the triangulation $K=(\tilde{K})^*$ as in Remark~\ref{oss:vl-Schauderhats}. Adopting the notation in that remark, If $k_{n-1}=0$, then $K=\tilde{K}$; else, we add the $0$-simplex $(x_{(n-1)n})$ and replace each $m$-simplex of the form $\sigma=(x_{u_0},\dots,x_{n-1},x_n)$ with the $m$-simplices $\tau=(x_{u_0},\dots,x_{n-1},x_{(n-1)n})$ and $\rho=(x_{u_0},\dots,x_{(n-1)n},x_n)$. The vl-Schauder hats of this new triangulation $K$ are $2k_0,\dots,2k_{n-1}$, together with the hats $h_i'=h_i-k_i$ for $i=0,\dots,n-1$ and $h_n^n$. Clearly, for all $i\neq j$, $2k_i\neq2k_j$ unless $k_i=k_j=0$: indeed, $2k_i$ attains its maximum at $x_{in}$, $2k_j$ attains its maximum at $x_{jn}$ and these two points are distinct because $x_i\neq x_j$.
 \end{proof}
\end{cl}

\begin{cl}
 If $\nu^+$ is a pc-valuation on $\P^+$ that assigns $1$ to each vl-Schauder hat, then $\nu^+(f)=\alpha^+(f)$ for all $f\in \P^+$.
 \begin{proof}
  We can write each $0\neq f\in \P^+$ as a sum $\sum_{i=0}^m a_ih_i$ (where $0<a_i\in\R$) of distinct vl-Schauder hats $h_0,\dots,h_m$ of a common triangulation $K$ that linearizes $f$. By Corollary~\ref{cor:vf}, we also have
  \begin{equation*}
 \nu^+(f)=\nu^+(\sum_{i=0}^m h_i)\qquad\text{and}\qquad\alpha^+(f)=\alpha^+(\sum_{i=0}^m h_i).
\end{equation*}
 We proceed by induction on  $m$.
If $m=0$, then, by Lemma~\ref{cl:alfavl-Schauder}, $$\nu^+(f)=\nu^+(h_1)=1=\alpha^+(h_1)=\alpha^+(f).$$

If $m>0$, by the induction hypotesis, for all $n<m$, $\nu^+\left(\sum_{j=0}^nl_j\right)=\alpha^+\left(\sum_{j=0}^nl_j\right)$ for distinct vl-Schauder hats $l_0,\dots,l_n$ of the same triangulation $H_n$ of $P$. Then, by Lemma~\ref{cl:andhat} and Corollary~\ref{cor:vf},
   \begin{align*}
    \nu^+(f)&=\nu^+\left(\sum_{i=0}^{m-1}h_i + h_m\right)=\\
    &=\nu^+\left(\sum_{i=0}^{m-1}h_i\right)+\nu^+(h_m)-\nu^+\left(\left(\sum_{i=0}^{m-1}h_i\right)\land h_m\right)=\\
    &=\nu^+\left(\sum_{i=0}^{m-1}h_i\right)+\nu^+(h_m)-\nu^+\left(\sum_{i=0}^{m-1}k_i\right)=\\
    &=\nu^+\left(\sum_{i=0}^{m-1}h_i\right)+\nu^+(h_m)-\nu^+\left(\sum_{i=0}^{m-1}2k_i\right)=\\
    &=\alpha^+\left(\sum_{i=0}^{m-1}h_i\right)+\alpha^+(h_m)-\alpha^+\left(\sum_{i=0}^{m-1}2k_i\right)=\\
    &=\alpha^+\left(\sum_{i=0}^{m-1}h_i\right)+\alpha^+(h_m)-\alpha^+\left(\sum_{i=0}^{m-1}k_i\right)=\\
    &=\alpha^+\left(\sum_{i=0}^{m-1}h_i\right)+\alpha^+(h_m)-\alpha^+\left(\left(\sum_{i=0}^{m-1}h_i\right)\land h_m\right)=\\
    &=\alpha^+\left(\sum_{i=0}^{m-1}h_i + h_m\right)=\\
    &=\alpha^+(f).
   \end{align*} 
 \end{proof}
\end{cl}

Now we can extend the pc-valuation $\alpha^+$ to a vl-valuation $\alpha$.
To complete the proof, we define $\alpha:\P\to\R$ as the functional such that, for all $f\in\P,$
$$\alpha(f)=\alpha^+(f^+)-\alpha^+(f^-).$$

Thanks to the uniqueness of the extension of a pc-valuation granted by Lemma~\ref{pcval}, $\alpha$ is the unique vl-valuation that assigns $1$ to each vl-Schauder hat of $\P$, and the Theorem  is proved.

\subsection*{Acknowledgement} The results presented here are part of the author's work towards a Ph.D. in Theoretical Computer Science at the Universit\`{a} degli Studi di Milano, Italy. The author wishes to thank Professor Vincenzo Marra for his help and 
suggestions, and for the supervision of the preparation of this paper.

\bibliographystyle{amsplain}

\providecommand{\bysame}{\leavevmode\hbox to3em{\hrulefill}\thinspace}
\providecommand{\MR}{\relax\ifhmode\unskip\space\fi MR }
\providecommand{\MRhref}[2]{%
  \href{http://www.ams.org/mathscinet-getitem?mr=#1}{#2}
}
\providecommand{\href}[2]{#2}

\end{document}